\documentclass[12pt,a4paper]{article}
\usepackage{amsmath,amsfonts,amssymb,amscd}

\usepackage{graphicx}
\usepackage{amsfonts,amssymb,amscd,amsmath}

\usepackage{indentfirst}

\usepackage{tikz}
\usetikzlibrary{matrix}

\def\H{\operatorname{H}}

\def\Ker{\operatorname{Ker}}
\def\Ad{\operatorname{Ad}}

\def\red{\operatorname{red}}

\def\GL{\operatorname{GL}}

\def\OSp{\operatorname{OSp}}
\def\Sp{\operatorname{Sp}}

\def\Im{\operatorname{Im}}

\newcounter{th}
\def\t{\refstepcounter{th}{\bf \noindent{Theorem} \arabic{th}. }}

\newcounter{prop}
\def\prop{\refstepcounter{prop}{\bf \noindent{Proposition} \arabic{prop}. }}

\newcounter{lem}
\def\lem{\refstepcounter{lem}{\bf \noindent{Lemma} \arabic{lem}. }}

\newcounter{de}
\def\de{\refstepcounter{de}{\bf \noindent{Definition} \arabic{de}. }}

\newcounter{ex}

\begin{document}

\begin{center}

{\LARGE{\bf Vector fields on $\mathfrak{osp}_{2m|2n}(\mathbb C)$- and $\pi\mathfrak{sp}_{n}(\mathbb C)$-flag supermanifolds}}\footnote[1]{Supported by  Universidade Federal de Minas Gerais, Brazil.}

\bigskip

\bigskip

{\bf E.G.Vishnyakova}\\[0.3cm]
\end{center}

\bigskip

\begin{abstract}
	The paper is devoted to a computation of the Lie superalgebras of holomorphic vector fields on  isotropic flag supermanifolds  of maximal type corresponding to the Lie superalgebras $\mathfrak{osp}_{2m|2n}(\mathbb C)$ and $\pi\mathfrak{sp}_{n}(\mathbb C)$. The result is that under some restrictions on the flag type any holomorphic vector field is fundamental with respect to the natural action of the Lie superalgebras $\mathfrak{osp}_{2m|2n}(\mathbb C)$ or $\pi\mathfrak{sp}_{n}(\mathbb C)$.

\end{abstract}

\bigskip

\section{Introduction} 
 Yu.I.~Manin \cite{Man} constructed four series of complex homogeneous supermanifolds that correspond to four series of  linear Lie superalgebras: $\mathfrak{gl}_{m|n}(\mathbb C)$, $\mathfrak{osp}_{m|2n}(\mathbb C)$, $\pi\mathfrak{sp}_{n}(\mathbb C)$ and $\mathfrak{q}_{n}(\mathbb C)$. (For definitions of these Lie superalgebras we refer to \cite{Kac}, see also Section $2$.) These supermanifolds are super-analogues of classical flag manifolds. In this paper we calculate the Lie superalgebras of global holomorphic vector fields on isotropic flag supermanifolds of maximal type that correspond to the Lie superalgebras $\mathfrak{osp}_{2m|2n}(\mathbb C)$ and $\pi\mathfrak{sp}_{n}(\mathbb C)$. We prove that under some restrictions on the flag type all such vector fields are fundamental with respect to the natural action of the corresponding Lie superalgebra.  We use induction and a similar result for super-Grassmannians that was obtained in \cite{oniosp,onipisp}. For isotropic flag supermanifolds and even for super-Grassmannians of non-maximal type an analogous result is not known so far. The Lie superalgebra of holomorphic vector  fields on super-Grassmannians corresponding to the Lie superalgebras $\mathfrak{gl}_{m|n}(\mathbb C)$ and $\mathfrak{q}_{n}(\mathbb C)$ were studied  in \cite{onigl,oniq}
and a similar question in the case of  flag supermanifolds  were studied in \cite{ViGL in JA,ViPi-sym}.

  The orthosymplectic Lie superalgebra $\mathfrak{osp}_{m|2n}(\mathbb C)$ is the linear Lie superalgebra that annihilates a non-degenerate  even symmetric bilinear form in $\mathbb C^{m|2n}$. The Lie superalgebra $\pi\mathfrak{sp}_{n}(\mathbb C)$ is the linear Lie superalgebra that annihilates a non-degenerate odd skew-symmetric bilinear form in $\mathbb C^{n|n}$. (For a detailed description of these Lie superalgebras see Section $2$.)

We denote by $\mathbf
 F_{k|l}$ the flag supermanifold of type $k|l$ in the vector superspace $\mathbb C^{m|n}$, see \cite{Man} and also \cite{ViGL in JA}. Here $k=(k_0,\ldots,k_r)$ and
 $l=(l_0,\ldots,l_r)$ such that
 \begin{equation}\label{eq conditions on k_i and l_j}
 	\begin{split}
 		0\le k_r\le\ldots\le k_0= m&,\quad 0\le l_r\ldots\le l_0= n,\\
 		\quad 0 <
 		k_r +l_r <\ldots& < k_0 + l_0 = m+n.
 	\end{split}
 \end{equation} 
This flag supermanifold corresponds to the Lie superalgebra  $\mathfrak{gl}_{m|n}(\mathbb C)$. The number $r$ is called the {\it length} of $\mathbf
 F_{k|l}$.  We denote by $\mathbf F^{e}_{k|l}$ the isotropic flag supermanifold in $\mathbb C^{m|2n}$ corresponding to $\mathfrak{osp}_{m|2n}(\mathbb C)$ and by $\mathbf
 F^{o}_{k|l}$ the isotropic flag supermanifold in $\mathbb C^{n|n}$ corresponding to $\pi\mathfrak{sp}_{n}(\mathbb C)$. Here the subscripts $e$ and $o$ in $\mathbf F^{e}_{k|l}$ and  $\mathbf F^{o}_{k|l}$ come from ``even'' and ``odd''.

  The idea of the proof is the following. For $r > 1$ the flag supermanifolds $\mathbf F^{e}_{k|l}$ and  $\mathbf F^{o}_{k|l}$ are the total spaces of holomorphic superbundles with the base spaces that are isomorphic to the isotropic super-Grassmannians and with the fibers that are isomorphic to the flag supermanifold $\mathbf F_{k'|l'}$ of length $r-1$. Here $k'=(k_1,\ldots,k_r)$ and
  $l'=(l_1,\ldots,l_r)$. Hence to obtain the result we can use induction and the results about Lie superalgebaras of holomorphic vector fields on super-Grassmannians \cite{oniosp,oniq} and on flag supermanifolds $\mathbf   F_{k'|l'}$ from \cite{ViGL in JA}.

  We set $\mathfrak{pgl}_{m|n}(\mathbb C):=\mathfrak{gl}_{m|n}(\mathbb C)/ \mathfrak
  z(\mathfrak{gl}_{m|n}(\mathbb C))$, where $\mathfrak
  z(\mathfrak{gl}_{m|n}(\mathbb C))$ is the center of $\mathfrak{gl}_{m|n}(\mathbb C)$. The main result of this paper was announced in \cite{Vivector}.  It is the following.

\medskip

\t\label{teor main osp} {\sl Let $r>1$.
	\begin{enumerate}
\item Assume that $m=k_1,\; n=l_1;\; (k_i,l_i)\ne (k_{i-1},0),\,(0,l_{i-1}),\;
i\geq 2,\; k_1\geq 1,\,l_1\geq 1$ and
$\mathfrak{v}(\mathbf{F}_{k'|l'})\simeq
\mathfrak{pgl}_{k_1|l_1}(\mathbb C)$. Then
$$
\mathfrak v(\mathbf F^{e}_{k|l})\simeq
\mathfrak {osp}_{2m|2n}(\mathbb C).
$$

\item Assume that $n=k_1+l_1;\; (k_i,l_i)\ne (k_{i-1},0),\,(0,l_{i-1}),\;
i\geq 2;,\; k_1\geq 3$, $l_1\geq 2$ and
$\mathfrak{v}(\mathbf{F}_{k'|l'})\simeq
\mathfrak{pgl}_{k_1|l_1}(\mathbb C)$. Then
$$
\mathfrak v(\mathbf F^{o}_{k|l})\simeq
\pi\!\mathfrak{sp}_{n}(\mathbb C).
$$
\end{enumerate}
	
}
\medskip

The Lie superalgebras of holomorphic vector fields on the flag super\-manifolds $\mathbf{F}_{k|l}$ and $\Pi\mathbf{F}_{k|k}$ corresponding to $\mathfrak{gl}_{m|n}(\mathbb C)$ and $\mathfrak{q}_{n}(\mathbb C)$, respectively, were calculated in \cite{ViGL in JA,ViPi-sym}. We obtained there the following result. 

\medskip

\t\label{teor main gl} {\sl {\bf A.} Assume that $r>1$ and that we have the following restrictions on the flag type: 
	\begin{align*}
	&(k_i,l_i)\ne (k_{i-1},0),\,(0,l_{i-1}),\,\, i\geq 2;\\
	&(k_{i-1},k_i|l_{i-1},l_i)\ne (1,0|l_{i-1},l_{i-1}-1),\, 	(1,1|l_{i-1},1),\,\, i\geq 1;\\
	&(k_{i-1},k_i|l_{i-1},l_i)\ne (k_{i-1},k_{i-1}-1|1,0),\, 	(k_{i-1},1|1,1),\,\, i\geq 1;\\
	&k|l\ne 	(0,\ldots,0| n,l_2,\ldots,l_r),\,\, k|l\ne 	( m,k_2,\ldots,k_r|0,\ldots,0).
	\end{align*}
	Then
	$$
	\mathfrak{v}(\mathbf{F}_{k|l})\simeq \mathfrak {pgl}_{m|n}(\mathbb
	C).
	$$
	If $k|l= (0,\ldots,0| n,l_2,\ldots,l_r)$ or $k|l= 	( m,k_2,\ldots,k_r|0,\ldots,0)$, then
	$$
	\mathfrak{v}(\mathbf{F}_{k|l})\simeq W_{mn}\subset\!\!\!\!\!\!+ (\bigwedge(\xi_1,\ldots,\xi_{mn})\otimes\mathfrak {pgl}_{n}(\mathbb
	C)),
	$$
	where $W_{mn}=\operatorname{Der}\bigwedge(\xi_1,\ldots,\xi_{mn})$.
	
	{\bf B.} Assume that $r > 1$, then for any $k$
$$
\mathfrak v(\Pi\mathbf{F}_{k|k})\simeq \mathfrak
{q}_{n}(\mathbb C)/\mathfrak{z}(\mathfrak{q}_{n}(\mathbb C)),
$$
 where $\mathfrak{z}(\mathfrak{q}_{n}(\mathbb C))$ is the center of $\mathfrak{q}_{n}(\mathbb C)$. 
}

\medskip

\section{Examples of Lie supergroups and Lie superalgebras}

A Lie supergroup is a group object in the category of supermanifolds. As in the classical Lie theory we can assign the Lie superalgebra to any Lie supergroup.  For more information about Lie supergroups see for example \cite{ViLieSupergroup}.  
Further we will need a description of some classical Lie supergroups and their Lie superalgebras. 

The general Lie supergroup $\operatorname{GL}_{m|n}(\mathbb C)$ is an open subsupermanifold in the superdomain
$$
\operatorname{Mat}_{m|n}(\mathbb C)= \left\{ \left(
\begin{array}{cc}
X_{11} & X_{12} \\
X_{21} & X_{22} \\
\end{array}
\right)\right\}.
$$
Here we consider the elements of the matrices $X_{11}\in
\operatorname{Mat}_m(\mathbb C)$ and $X_{22}\in
\operatorname{Mat}_n(\mathbb C)$ as even coordinates of the superdomain $\operatorname{Mat}_{m|n}(\mathbb C)$ and the elements of the matrices $X_{12}$, $X_{21}$ as odd ones. The subsupermanifold $\operatorname{GL}_{m|n}(\mathbb C)$ is defined by the following
equations 
$$
\det X_{11}\ne 0\quad \text{and} \quad
\det X_{22}\ne 0.
$$ 
The multiplication in the Lie supergroup $\operatorname{GL}_{m|n}(\mathbb C)$ is given by the usual matrix multiplication.

The Lie superalgebra $\mathfrak{gl}_{m|n}(\mathbb
C)$ of $\operatorname{GL}_{m|n}(\mathbb C)$ has the following form:
$$
\mathfrak{gl}_{m|n}(\mathbb C)= \left\{ \left(
\begin{array}{cc}
A_{11} & A_{12} \\
A_{21} & A_{22} \\
\end{array}
\right)\right\},
$$
where $A_{11}$, $A_{12}$, $A_{21}$ and $A_{22}$ are matrices of complex numbers of size
$m\times m$, $m\times n$, $n\times m$ and $n\times n$, respectively. The even part 
$\mathfrak{gl}_{m|n}(\mathbb C)_{\bar 0}$ of $\mathfrak{gl}_{m|n}(\mathbb
C)$ is determined by the equations 
$A_{12}=0$, $A_{21}=0$, and the odd part
$\mathfrak{gl}_{m|n}(\mathbb C)_{\bar 1}$ is given by $A_{11}=0$,
$A_{22}=0$. The multiplication in this Lie superalgebra has the form: 
$$
[X,Y]=XY-(-1)^{p(X)p(Y)}YX,
$$
where $X,\,Y$ are homogeneous elements in $\mathfrak{gl}_{m|n}(\mathbb C)$ and 
$p(Z)$ is the parity of $Z$. The center $\mathfrak
z(\mathfrak{gl}_{m|n}(\mathbb C))$ of $\mathfrak{gl}_{m|n}(\mathbb C)$ contains all matrices $\alpha E_{m+n}$, where $\alpha\in \mathbb C$ and $E_{m+n}$ is the identity matrix of size $m+n$. By definition we put $\mathfrak{pgl}_{m|n}(\mathbb C):=\mathfrak{gl}_{m|n}(\mathbb C)/ \mathfrak
z(\mathfrak{gl}_{m|n}(\mathbb C))$. 

\medskip

Consider the following two classical Lie subsuperalgebras in $\mathfrak{gl}_{m|n}(\mathbb C)$.

\medskip

\noindent{\bf (1)} The Lie superalgebra $ \mathfrak {osp}_{m|2n}(\mathbb C)$ is a Lie subsuperalgebra in $\mathfrak{gl}_{m|2n}(\mathbb C)$ that annihilates a non-degenerate even symmetric bilinear form $\beta$ in $\mathbb
C^{m|2n}$. The matrix $\Gamma$ of $\beta$ in the standard basis in  $\mathbb C^{m|2n}$ for even and odd  $m$ is given respectively by
\begin{equation}
\label{matr_grama_osp}
\Gamma=\left(
\begin{array}{cccc}
0 & E_s & 0 & 0 \\
E_s & 0 & 0 & 0 \\
0 & 0 & 0 & E_n \\
0 & 0 & -E_n & 0 \\
\end{array}
\right), \quad \Gamma=\left(
\begin{array}{ccccc}
0&E_s & 0 & 0 & 0 \\
E_s& 0 & 0 & 0 & 0 \\
0&0 & 1 & 0 & 0 \\
0&0 & 0 & 0 & E_n \\
0&0 & 0 & -E_n & 0 \\
\end{array}
\right).
\end{equation}
Here
$m=2s$ or $m=2s+1$. Explicitly we have
\begin{equation}
\label{matr_osp_chet} \mathfrak {osp}_{2s|2n}(\mathbb C)=\left\{
\left(
\begin{array}{cccc}
A_{11} & A_{12} & C_{11} & C_{12} \\
A_{21} & -A^T_{11} & C_{21} & C_{22} \\
-C_{22}^T & -C_{12}^T & B_{11} & B_{12} \\
C_{21}^T & C_{11}^T & B_{21} & -B_{11}^T \\
\end{array}
\right),\,\,
\begin{array}{l}
A_{21}^T=-A_{21},\\ A_{12}^T=-A_{12}, \\
B_{12}^T=B_{12}, \\B_{21}^T=B_{21}
\end{array}
\right\},
\end{equation}
and
\begin{equation}
\label{matr_osp_nechet} \mathfrak {osp}_{2s+1|2n}(\mathbb C)=
\left\{ \left(
\begin{array}{ccccc}
A_{11}&A_{12}&G_1&C_{11}& C_{12}\\
A_{21} &-A_{11}^T & G_{2} & C_{21} & C_{22} \\
-G_2^T &   -G_1^T & 0 & G_3 & G_4 \\
-C_{22}^T & -C_{12}^T & -G^T_{4} & B_{11} & B_{12} \\
C_{21}^T &   C_{11}^T & G_{3}^T & B_{21} & -B_{11}^T \\
\end{array}
\right),
\begin{array}{l}
A_{21}^T=-A_{21},\\A_{12}^T=-A_{12},\\
B_{12}^T=B_{12}, \\ B_{21}^T=B_{21}
\end{array}
\right\}.
\end{equation}
Here $A_{11}$, $B_{11}$ are square matrices of size $s$ and $n$, respectively. The center $\mathfrak z(\mathfrak {osp}_{m|2n}(\mathbb
C))$ of $\mathfrak {osp}_{m|2n}(\mathbb
C)$ is trivial. The corresponding connected Lie supergroup  we will denote by $\operatorname{OSp}_{m|2n}(\mathbb C)$. This is a subsupermanifold in $\operatorname{GL}_{m|2n}(\mathbb C)$ that is given by the following equation:
$$
\left(
\begin{array}{cc}
X_{11} & X_{12} \\
X_{21} & X_{22} \\
\end{array}
\right)^{ST}\Gamma \left(
\begin{array}{cc}
X_{11} & X_{12} \\
X_{21} & X_{22} \\
\end{array}
\right)=\Gamma,
$$
where
$$
\left(
\begin{array}{cc}
X_{11} & X_{12} \\
X_{21} & X_{22} \\
\end{array}
\right)^{ST}=\left(
\begin{array}{cc}
X_{11}^T & X_{21}^T \\
-X_{12}^T & X_{22}^T \\
\end{array}
\right).
$$
and $T$ is the usual transposition. 

\medskip

{\bf (2)} The Lie superalgebra $\pi\mathfrak{sp}_{n}(\mathbb C)\subset \mathfrak
{gl}_{n|n}(\mathbb C)$ is a Lie subsuperalgebra in $\mathfrak{gl}_{n|n}(\mathbb C)$ that annihilate a non-degenerate odd skew-symmetric bilinear form
 $\gamma$ in $\mathbb C^{n|n}$. The matrix $\Upsilon$ of $\gamma$ in the standard basis in $\mathbb C^{n|n}$ has the following form:
\begin{equation}
\label{matr_grama_pisp}
\Upsilon=\left(
\begin{array}{cc}
0 & E_n \\
E_n & 0 \\
\end{array}
\right).
\end{equation}
Then we have
\begin{equation}
\label{matr_pisp} \pi\mathfrak{sp}_{n}(\mathbb C)=\left\{\left(
\begin{array}{cc}
A & B \\
C & -A^T \\
\end{array}
\right), \quad B^T=B,\,\, C^T=-C\right\},
\end{equation}
where $A$, $B$, $C$ are square matrices of size $n$. The center 
$\mathfrak z(\pi\mathfrak{sp}_{n}(\mathbb C))$ of $\pi\mathfrak{sp}_{n}(\mathbb C)$ is trivial. We will denote the corresponding connected Lie supergroup   by $\Pi\!\operatorname{Sp}_{n}(\mathbb C)$. This is a subsupermanifold in $\operatorname{GL}_{n|n}(\mathbb C)$ that is given by the following equation:
$$
\left(
\begin{array}{cc}
X_{11}& X_{12} \\
X_{21} & X_{22} \\
\end{array}
\right)^{ST}\Upsilon \left(
\begin{array}{cc}
X_{11} & X_{12} \\
X_{21} & X_{22} \\
\end{array}
\right)=\Upsilon.
$$

\section{Flag supermanifolds}
 
An introduction to the theory of supermanifolds can be found in \cite{BL,Fioresi,Leites,Man}.  Throughout this paper we will be interested in the complex-analytic version of the theory. A {\it complex-analytic supermanifold} of dimension $p|q$ is a $\mathbb{Z}_2$-graded ringed space that
is locally isomorphic to a complex-analytic superdomain of dimension $p|q$, this is to a ringed space of the form $\mathcal U = (\mathcal U_0, \mathcal{F}_{\mathcal U_0} \otimes_{\mathbb C} \bigwedge(q) )$. Here $\mathcal{F}_{\mathcal U_0}$ is the sheaf of holomorphic functions on an open set $\mathcal U_0\subset \mathbb{C}^p$ and
$ \bigwedge(q)$ is the Grassmann algebra with $q$ generators. We will denote a supermanifold by $\mathcal{M} = (\mathcal{M}_0,{\mathcal O}_{\mathcal{M}})$, where $\mathcal{M}_0$ is the underlying complex-analytic manifold and ${\mathcal O}_{\mathcal{M}}$ is the structure sheaf of $\mathcal{M}$. A {\it morphism of complex-analytic} supermanifolds
$\mathcal M$ to  $\mathcal N$ is a pair $f = (f_0,f^*)$, where $f_0: \mathcal M_0\to \mathcal N_0$ is a holomorphic map and $f^*: \mathcal O_{\mathcal N}\to (f_0)_*(\mathcal
O_{\mathcal M})$ is a homomorphism of sheaves of superalgebras. 
 We denote by $\mathcal T = \mathcal Der\,(\mathcal O_{\mathcal M})$  the sheaf of holomorphic vector fields on  $\mathcal M$. The sheaf $\mathcal T$ is a sheaf of Lie superalgebras with respect to the following multiplication 
$$
[X,Y] = XY
-(-1)^{p(X)p(Y)}YX.
$$
The global sections of $\mathcal T$
are called {\it holomorphic vector fields} on $\mathcal M$. We will denote the Lie superalgebra of global holomotphic vector fields  by $\mathfrak v(\mathcal M)$. It is known that the Lie superalgebra $\mathfrak v(\mathcal M)$ is finite dimensional if $\mathcal M_0$ is compact. (A proof of this fact can be found for instance in \cite{ViGL in JA}.)

\medskip

\noindent{\bf $\mathfrak {gl}_{m|n}(\mathbb C)$-flag supermanifolds.} We fix two sets of non-negative integers $k=(k_0,\ldots,k_r)$ and $l=(l_0,\ldots,l_r)$ such that (\ref{eq conditions on k_i and l_j}) holds. Let us recall the definition of a flag supermanifold $\mathbf F_{k|l}$ in the superspace $\mathbb C^{m|n}$, see \cite{Man} and also \cite{ViGL in JA}. We use definitions and notations from \cite{ViGL in JA}. The only difference is that in this paper for convenience we denote by $\mathbf F_{k|l}$ the flag supermanifold of type $k|l$, where  $k=(k_0,\ldots,k_r)\quad  \text{and}$ $l=(l_0,\ldots,l_r)$
and we assume that $k_0=m$ and $l_0=n$. In \cite{ViGL in JA} this supermanifold was denoted by $\mathbf F_{k|l}^{m|n}$. 
The underlying space of the supermanifold
$\mathbf F_{k|l}$ is the product  $\mathbf
F_{k}\times \mathbf F_{l}$ of two flag manifolds of types $k$ and $l$, respectively.

Let us describe an atlas on $\mathbf F_{k|l}$, see \cite{ViGL in JA} for details. We fix two subsets $I_{s\bar
0}\subset\{1,\ldots,k_{s-1}\}$  and $I_{s\bar 1}\subset\{1,\ldots,l_{s-1}\}$
  such that $|I_{s\bar 0}| = k_s,$ and $|I_{s\bar 1}| = l_s$,  where $s = 1,\ldots,r$, and we put $I_s:=(I_{s\bar 0},I_{s\bar 1})$ and $I := (I_1,\ldots,I_r)$. To any such $I_s$  we assign the following $(k_{s-1} + l_{s-1})\times (k_s + l_s)$-matrix 
\begin{equation}\label{eq local chart general}
Z_{I_s}
=
\left(
\begin{array}{cc}
X_s & \Xi_s\\
\H_s & Y_s \end{array} \right), \ \ s=1,\dots,r.
\end{equation}
 We assume that the matrices $X_s=(x^s_{ij})$ and $Y_s=(y^s_{ij})$ in (\ref{eq local chart general}) have size $(k_{s-1}\times k_s)$ and $(l_{s-1}\times l_s)$, respectively,  and that $Z_{I_s}$ contains the identity submatrix $E_{k_s+l_s}$ of size $(k_s+l_s)\times (k_s+l_s)$ in the lines with numbers $i\in I_{s\bar 0}$ and $k_{s-1} + i,\; i\in I_{s\bar 1}$.   For example, if 
 $	I_{s\bar 0}=\{k_{s-1}-k_s+1,\ldots, k_{s-1}\},$ $I_{s\bar 1}=\{l_{s-1}-l_s+1,\ldots, l_{s-1}\}
$, then
 the matrix $Z_{I_s}$ has the following form:
 $$
 Z_{I_s} =\left(
 \begin{array}{cc}
 X_s&\Xi_s\\
 E_{k_s}&0\\
 \H_s&Y_s\\0&E_{l_s}\end{array} \right).
 $$
 Here we denote by $E_{q}$ the identity square matrix of size $q$. To simplify notations we use the same letters $X_s$, $Y_s$, $\Xi_s$ and $\H_s$ as in (\ref{eq local chart general}).

 The matrices (\ref{eq local chart general}) determine the superdomain $\mathcal U_I$ with even coordinates $x^s_{ij}$ and $y^s_{ij}$, and odd coordinates $\xi^s_{ij}$ and $\eta^s_{ij}$. The transition functions between two  superdomains $\mathcal U_I$ and $\mathcal U_J$, where $I = (I_s)$ and $J = (J_s)$, are defined in the following way:
\begin{equation}\label{eq transition functions}
Z_{J_1} = Z_{I_1}C_{I_1J_1}^{-1}, \quad Z_{J_s} =
C_{I_{s-1}J_{s-1}}Z_{I_s}C_{I_sJ_s}^{-1},\quad s\ge 2. 
\end{equation}
The matrix $C_{I_1J_1}$ is an invertible submatrix in $Z_{I_1}$ that consists of the lines with numbers $i\in J_{1\bar 0}$ and $k_{0} + i,$ where $i\in J_{1\bar 1}$, and  $C_{I_sJ_s}$ is the invertible submatrix in $C_{I_{s-1}J_{s-1}}Z_{I_s}$ that consists of the lines with numbers $i\in J_{s\bar 0}$ and $k_{s-1} + i,$ where $i\in J_{s\bar 1}$, see \cite{ViGL in JA} for details.  Now the atlas on $\mathbf F_{k|l}$ is described. The supermanifold $\mathbf F_{k|l}$ is called the {\it flag supermanifold  of type $k|l$}. In case $r = 1$ this supermanifold is called the {\it super-Grass\-mannian} and sometimes it is denoted in the literature by $\mathbf {Gr}_{m|n,k|l}$.

Recall that we denote by $\mathfrak {gl}_{m|n}(\mathbb C)$ the general Lie superalgebra of the superspace $\mathbb C^{m|n}$ and by $\GL_{m|n}(\mathbb C)$ the connected Lie supergroup of the Lie superalgebra $\mathfrak {gl}_{m|n}(\mathbb C)$.
 In \cite{Man} an action of  $\GL_{m|n}(\mathbb C)$ on
 $\mathbf {F}_{k|l}$ is defined. Let us recall this definition in our notations and in our atlas. Let 
 $$
 L= 
 \left(
 \begin{array}{cc}
 L_{11} & L_{12} \\
 L_{21} & L_{22} \\
 \end{array}
 \right)
 $$
 be a coordinate matrix of the Lie supergroup $\GL_{m|n}(\mathbb C)$. Then the action of  $\GL_{m|n}(\mathbb C)$ on 
 $\mathbf {F}_{k|l}$ in our coordinates is given by the following formulas:
\begin{equation}\label{eq action of Q}
\begin{aligned} 
&(L,(Z_{I_1},\ldots,Z_{I_r})) \longmapsto
(\tilde Z_{J_1},\ldots,\tilde Z_{J_r}), 
\end{aligned} 
\end{equation}
where $\tilde Z_{J_1} =
LZ_{I_1}C_1^{-1}$ and $\tilde Z_{J_s} = C_{s-1}Z_{I_s}C_s^{-1}$. 
Here $C_1$ is the invertible submatrix in $LZ_{I_1}$ that consists of the lines with numbers $i\in J_{1\bar 0}$ and $m+i$, where  $i\in J_{1\bar 1}$, and $C_s,$ where $s\ge 2$, is the invertible submatrix in $C_{s-1}Z_{I_s}$ that consists of the lines with numbers $i\in J_{s\bar 0}$ and $k_{s-1}+i$, where
$i\in J_{\bar 1s}$. This Lie supergroup action induces the following Lie superalgebra homomorphism 
$$
\mu:\mathfrak {gl}_{m|n}(\mathbb C)\to\mathfrak v(\mathbf {F}_{k|l}).
$$

\medskip

\noindent{\bf  $\mathfrak
	{osp}_{m|2n}(\mathbb C)$- and $\pi\mathfrak{sp}_{n}(\mathbb C)$-flag supermanifolds.} We can assign to the Lie superalgebras  $\mathfrak
{osp}_{m|2n}(\mathbb C)$ and $\pi\mathfrak{sp}_{n}(\mathbb C)$ isotropic flag supermanifolds 
$\mathbf F_{k|l}^{e}$ and $\mathbf F_{k|l}^{o}$, respectively. 
The underlying spaces of these supermanifolds are the manifolds of isotropic flags with respect to the form  $\beta$ and
$\gamma$, see Section $2$. In case $r=1$ these supermanifolds are called  {\it isotropic super-Grassmannians}. These supermanifolds were studied in \cite{oniosp,oniq}.  Let us describe the supermanifolds $\mathbf F_{k|l}^{e}$ and $\mathbf F_{k|l}^{o}$ using charts and local coordinates.

 The supermanifold $\mathbf F_{k|l}^{e}$, where $l_0$ is even, is a subsupermanifold of  $\mathbf F_{k|l}$ that is given in local coordinates $(\ref{eq local chart general})$ by the following equation:
\begin{equation}
\begin{split}
\left(
\begin{array}{cc}
X_1 & \Xi_1\\
\H_1 & Y_1 \end{array} \right)^{ST} \Gamma \left(
\begin{array}{cc}
X_1 & \Xi_1\\
\H_1 & Y_1 \end{array} \right)=0,
\end{split}
\end{equation}
where $\Gamma$ is as in (\ref{matr_grama_osp}) and $ST$ is the supertransposition, see Section $2$.  In case $m=2k_1$ or $m=2k_1+1$ and $n=l_1$, we say that the supermanifold $\mathbf F_{k|l}^{e}$ has {\it maximal type}. The subsupermanifold $\mathbf F_{k|l}^{o}$ of the supermanifold $\mathbf
F_{k|l}$ is given in the local coordinates $(\ref{eq local chart general})$ by the following equation:
\begin{equation}
\begin{split}
\left(
\begin{array}{cc}
X_1 & \Xi_1\\
\H_1 & Y_1 \end{array} \right)^{ST} \Upsilon \left(
\begin{array}{cc}
X_1 & \Xi_1\\
\H_1 & Y_1 \end{array} \right)=0,
\end{split}
\end{equation}
where $\Upsilon$ is as (\ref{matr_grama_pisp}). In case $n=k_1+l_1$, we say that  the supermanifold $\mathbf F_{k|l}^{o}$ has {\it maximal type}.

There are transitive actions $\mu_{e}$ and $\mu_{o}$ of the Lie supergroups  $\OSp_{m|2n}(\mathbb C)$ and $\Pi\!\Sp_n(\mathbb C)$ on $\mathbf F_{k|l}^{e}$ or $\mathbf F_{k|l}^{o}$, respectively. It is given by  $(\ref{eq action of Q})$, if we replace  $L$ by a coordinate matrix of the Lie supergroup $\OSp_{m|2n}(\mathbb C)$ or $\Pi\!\Sp_n(\mathbb C)$, respectively. This action induces the homomorphism of Lie superalgebras
$$
\mu_{e}:\mathfrak
{osp}_{m|2n}(\mathbb C)\to\mathfrak v(\mathbf F_{k|l}^{e})\quad \text{and} \quad \mu_{o}:\pi\mathfrak{sp}_{n}(\mathbb C)\to\mathfrak v(\mathbf F_{k|l}^{o}).
$$

\section{Vector fields on a superbundle}

For computation of the Lie superalgebra of holomorphic vector fields on isotropic flag supermanifolds we will use the following fact. For $r> 1$ the isotropic flag supermanifold $\mathbf {F}_{k|l}^e$ is a superbundle  with the base space that is isomorphic to the isotropic super-Grassmannian $\mathbf {F}_{k_0,k_1|l_0,l_1}^e$  and with the fiber that is isomorphic to $\mathbf F_{k'|l'}$. Similarly, for $r> 1$ the isotropic flag supermanifold $\mathbf {F}_{k|l}^o$ is a superbundle  with the base spaces  $\mathbf {F}_{k_0,k_1|l_0,l_1}^o$ and the fiber $\mathbf F_{k'|l'}$. (A similar statement holds for $\mathbf {F}_{k|l}$, see \cite{ViGL in JA} for details.)
 In local coordinates that we introduced the bundle projection say $\pi$ is given by $(Z_1,Z_2,\ldots Z_n) \longmapsto (Z_1).
$ Moreover, from Formulas (\ref{eq action of Q}) we can deduce that the projection $\pi$ is equivariant with respect to the action of the Lie supergroups $\OSp_{m|2n}(\mathbb C)$ and $\Pi\!\Sp_n(\mathbb C)$ on $\mathbf {F}_{k|l}^e$ and $\mathbf {F}_{k|l}^o$, respectively.

We will need some facts about vector fields on superbundles. 

\medskip
\de \label{de projectable vector fields} Let $\pi=(\pi_0,\pi^*): \mathcal M\to \mathcal N$ be a morphism of holomorphic supermanifolds. A vector field $v\in\mathfrak
v(\mathcal M)$ is called {\it projectible} with respect to $\pi$,
if there exists a vector field $v_1\in\mathfrak v(\mathcal N)$
such that $
\pi^*(v_1(f))=v(\pi^*(f))$ for all $f\in \mathcal O_{\mathcal N}$.
In this case we say that {\it $v$ is projected to} $v_1$.

\medskip

Denote the Lie superalgebra of projectible vector fields by  $\overline{\mathfrak v}(\mathcal M)$. Assume that
$\pi=(\pi_0,\pi^*): \mathcal M\to \mathcal N$ is the projection of a holomorphic superbundle. Then the homomorphism
$\pi^*: \mathcal O_{\mathcal N}\to \pi_*(\mathcal O_{\mathcal M})$ is injective and any projectible vector field $v$ is projected into a unique vector field
 $v_1 = \mathcal{P}(v)$. Hence we have the following homomorphism of Lie superalgebras 
 $$
 \mathcal{P}:
\overline{\mathfrak v}(\mathcal M)\to \mathfrak v(\mathcal N),\quad v \mapsto v_1.
$$
 A vector field $v\in\mathfrak v(\mathcal M)$ is called {\it vertical}, if $\mathcal{P}(v)=0$. The Lie subalgebra   $\Ker\mathcal{P}$ is an ideal in
 $\overline{\mathfrak v}(\mathcal M)$.

If $\mathcal S$ is a supermanifold, then the global sections $\mathcal O_{\mathcal S}(\mathcal S_0)$ of the structure sheaf $\mathcal O_{\mathcal S}$ are also called {\it holomorphic functions} on $\mathcal S$. 

\medskip

\prop\cite{Bash}\label{prop Bash} {\sl Let  
	$p: \mathcal M 	\to \mathcal B$ be the projection of a superbundle with fiber $\mathcal S = (\mathcal S_0, \mathcal O_{\mathcal S})$. Assume that $\mathcal O_{\mathcal S}(\mathcal S_0) = \mathbb C$, then any global holomorphic vector field on $\mathcal M$ is projectible.
	}

\medskip

Let $\pi: \mathcal M\to
\mathcal B$ be a superbundle with the fiber $\mathcal S$. We denote by $\mathcal W$ the sheaf of vertical holomorphic vector fields. It is a sheaf of the base space $\mathcal B_0$ of $\mathcal B$. More precisely, to any open set $U\subset \mathcal B_0$ we assign
 the set of all vertical vector fields on the supermanifold $(\pi_0^{-1}(U),\mathcal
O_{\mathcal M})$. 

\medskip
\prop\cite{Vi}\label{prop W is localy free}  {\sl Assume that the base space $\mathcal S_0$ of the fiber $\mathcal S$ is compact. Then $\mathcal W$ is a locally free sheaf of $\mathcal O_{\mathcal B}$-modules and
 $\dim\mathcal
W= \dim\mathfrak v(\mathcal S)$.}

\medskip

By definition we have $\mathcal W(\mathcal B_0)=\Ker\mathcal{P}$. In \cite{ViGL in JA} we described the corresponding to $\mathcal W$ graded sheaf $\tilde {\mathcal W}$. It is defined in the following way:
\begin{equation}\label{eq tilda W def}
\tilde {\mathcal W}=\bigoplus_{p\ge 0}\tilde{\mathcal W}_{p}, \quad
\text{where}\quad
\tilde{\mathcal W}_{p}=\mathcal W_{(p)}/\mathcal
W_{(p+1)}.
\end{equation}
Here $\mathcal W_{(p)}=
\mathcal J^p\mathcal W$ and $\mathcal J$ is the sheaf of ideals in ${\mathcal O}_{\mathcal B}$ generated by odd elements. 
Clearly, $\tilde {\mathcal W}$ is the $\mathbb Z$-graded 
sheaf of $\mathcal F_{\mathcal B_0}$-modules, where  $\mathcal F_{\mathcal B_0}$ is the structure sheaf of the underlying space $\mathcal B_0$. 
By Proposition \ref{prop W is localy free} we get the following result.

\medskip
\prop\label{prop tilde W_0 is localy free}  {\sl Assume that the base space $\mathcal S_0$ of the fiber $\mathcal S$ is compact. Then $\tilde{\mathcal W}_0$ is a locally free sheaf of $\mathcal F_{\mathcal B_0}$-modules and any fiber of the corresponding vector bundle is isomorphic to $\mathfrak v(\mathcal S)$.
	
}
\medskip

Denote by $\mathbf W_0$ the vector bundle corresponding to the locally free sheaf  $\tilde{\mathcal W}_0$. 
To calculate the Lie superalgebra of holomorphic vector fields on isotropic flag supermanifolds we will use Proposition \ref{prop Bash} and the following  result. 

\medskip
\t\cite{Viholom}\label{teor constant functions} {\sl Consider the flag supermanifold $\mathcal M=\mathbf{F}_{k|l}$.  Assume that 
	\begin{equation}\label{eq_condition no functions}
		\begin{split}
			(k|l)&\ne 
			(m,\ldots, m, k_{s+2},\ldots, k_r| l_1,\ldots, l_s,0,\ldots,0),\\
			(k|l)&\ne (k_1,\ldots, k_s,0,\ldots,0 |n,\ldots, n, l_{s+2},\ldots,
			l_r),
		\end{split}
	\end{equation}
	for any $s\geq 0$. Then $\mathcal O_{\mathcal M}(\mathcal M_0) = 	\mathbb{C}.$ 	Otherwise
	 $\mathcal O_{\mathcal M}(\mathcal M_0) = \bigwedge(mn)$, where $\bigwedge(mn)$ is the Grassmann algebra with $mn$ generators.$\Box$ 
}

\medskip

\section{The Borel-Weil-Bott Theorem}

Further we will use the Borel-Weil-Bott Theorem to compute   the vector space of  global sections of locally free sheaves. 
 Details about the Borel-Weil-Bott Theorem can be found for example in \cite{ADima}. Recall that this theorem is used to  compute cohomology group with values in a homogeneous holomorphic bundle over a product of classical or isotropic flag manifolds.   For completeness we formulate this theorem here adapting to our notations.  First of all we need to describe the underlying space of an isotropic super-Grassmannian of maximal type. 

\medskip

\noindent{\bf $\mathfrak {osp}_{m|2n}(\mathbb C)$-super-Grassmannians of maximal type.} Consider the isotropic super-Grassmanian of maximal type $\mathbf {F}_{m,s|2n,n}^e$, where $m=2s$ or $2s+1$. Denote by  $G^e_{2s} \simeq \operatorname{SO}_{2s}(\mathbb C)\times \operatorname{Sp}_{2n}(\mathbb
C)$ the underlying space of the Lie supergroup $\operatorname{OSp}_{2s|2n}(\mathbb C)$ and by $P^e_{2s}$ the parabolic subgroup in $G^e_{2s}$ that contains all matrices in the following form: 
\begin{equation}\label{eq parabolic OSP}
\left(
	\begin{array}{cccc} A_1 & 0& 0& 0\\ C_1& (A_1^T)^{-1}& 0& 0\\
		0& 0&A_2 &0\\
		0& 0& C_2& (A_2^T)^{-1}
	\end{array} \right),
\end{equation}
where $A_1\in \GL_{s}(\mathbb C)$ and $A_2\in \GL_{n}(\mathbb C)$. We also denote by $R^e_{2s}$ the reductive part of $P^e_{2s}$. Clearly, $R^e_{2s} \simeq \GL_{s}(\mathbb C)\times \GL_{n}(\mathbb C)$. The underlying manifold of $\mathbf {F}_{2s,s|2n,n}^e$ is isomorphic to $G^e_{2s}/P^e_{2s}$. We see that it is a product of two isotropic Grassmannians.

Further, the underlying manifold of $\mathbf {F}_{2s+1,s|2n,n}^e$ is isomorphic to the homogeneous space $G^e_{2s+1}/P^e_{2s+1}$,
 where $G^e_{2s+1} \simeq \operatorname{SO}_{2s+1}(\mathbb C)\times \operatorname{Sp}_{2n}(\mathbb
 C)$ is the underlying space of 
$\operatorname{OSp}_{2s+1|2n}(\mathbb C)_{\red}$ and $P^e_{2s+1}$ is the parabolic subgroup in $G^e_{2s+1}$ that contains all matrices in the following form:
\begin{equation}
\label{P_osp_m-1}
\left(
\begin{array}{ccccc}
A_1 & 0& 0& 0& 0\\
C_1& (A_1^T)^{-1}& G& 0& 0\\
H&0&1& 0& 0\\
0&0&0& A_2 & 0\\
0&0&0& C_2& (A_2^T)^{-1}
\end{array}
\right).
\end{equation}
Here 
$A_1\in \GL_{s}(\mathbb C)$ and $A_2\in \GL_{n}(\mathbb C)$.
The reductive part $R^e_{2s+1}$ of $P^e_{2s+1}$ has the form
$R^e_{2s+1}\simeq
\GL_{s}(\mathbb C)\times \GL_{n}(\mathbb C)$.

\medskip

\noindent{\bf $\pi\mathfrak{sp}_{n}(\mathbb C)$-super-Grassmannians of maximal type.} The underlying mani\-fold of the isotropic super-Grassmannian of maximal type $\mathbf {F}_{n,s|n,t}^o$, where $s+t=n$, is the usual Grassmannian $\mathbf {F}_{n,s}$, see \cite{Man,onipisp}. It is isomorphic to $G^o/P^o$, where $G^o\simeq \GL_n(\mathbb C)$ is the underlying space of $\Pi\!\Sp_{n}(\mathbb C)$  and
$P^o$ is the parabolic subgroup in $G^o$ that contains all matrices in the following form:
\begin{equation}
\label{P_pisp}
\left(
\begin{array}{cccc} A_1 & 0& 0& 0\\ C_1& B_1& 0& 0\\
0& 0& (A_1^T)^{-1}&D_1\\
0& 0&0& (B_1^T)^{-1} \end{array}
\right),
\end{equation}
where $A_1\in \GL_{t}(\mathbb C)$ and $B_1\in \GL_{s}(\mathbb C)$.
The reductive part of $P^o$ has the form $R^o \simeq \GL_{t}(\mathbb C)\times \GL_{s}(\mathbb C)$.

\medskip

To use the Borel-Weil-Bott Theorem we need to fix Cartan subalgebras and  root systems. 

\medskip

\noindent{\bf Cartan subalgebras and root systems.} 
In the Lie algebra  $\mathfrak {osp}_{2s|2n}(\mathbb C)_{\bar 0}=\mathfrak{so}_{2s}(\mathbb
C)\oplus\mathfrak {sp}_{2n}(\mathbb C)$
we fix the following Cartan subalgebra 
$$
\mathfrak t(\mathfrak {osp}_{2s|2n}(\mathbb C)_{\bar 0})=\mathfrak
t_0(\mathfrak {so}_{2s}(\mathbb C))\oplus\mathfrak t_1(\mathfrak
{sp}_{2n}(\mathbb C)),
$$
where
$$
\begin{array}{l}
\mathfrak t_0(\mathfrak {so}_{2s}(\mathbb C))=
\{\operatorname{diag}(\mu_1,\dots,\mu_s, -\mu_1,\dots,-\mu_s)\},\\
\mathfrak t_1(\mathfrak {sp}_{2n}(\mathbb C))=
\{\operatorname{diag}(\lambda_1,\dots,\lambda_n,-\lambda_1,
\dots,-\lambda_n)\}.
\end{array}
$$
We fix the following system of positive roots:
$$
\Delta^+(\mathfrak {osp}_{2s|2n}(\mathbb C)_{\bar
	0})=\Delta^+_1(\mathfrak {so}_{2s}(\mathbb C)) \cup
\Delta^+_2(\mathfrak {sp}_{2n}(\mathbb C)),
$$
where
$$
\begin{array}{l}
\Delta^+_1(\mathfrak {so}_{2s}(\mathbb C))=\{\mu_i-\mu_j,\,\,
\mu_i+\mu_j, \,\,i<j\},
\\
\Delta^+_2(\mathfrak {sp}_{2n}(\mathbb C))=\{\lambda_p-\lambda_q,\,\,
\,\,p<q, \,\, \lambda_p+\lambda_q,\,\, \,\,p\leq q \},
\end{array}
$$
and the following system of simple roots 
$$
\Phi(\mathfrak {osp}_{2s|2n}(\mathbb
C)_{\bar 0})=\Phi_1(\mathfrak {so}_{2s}(\mathbb C))\cup
\Phi_2(\mathfrak {sp}_{2n}(\mathbb C)),
$$
 where
 \begin{align*}
& \Phi_1(\mathfrak {so}_{2s}(\mathbb C))= \{\alpha_1,..., \alpha_{s}\},
 \,\,\, \alpha_i=\mu_i-\mu_{i+1},
 \,\,i=1,\ldots, s-1, \,\,\alpha_s=\mu_{s-1}+\mu_{s},\\
& \Phi_2(\mathfrak {sp}_{2n}(\mathbb C))= \{\beta_1,..., \beta _{n}\},
 \,\,\, \beta_j=\lambda_j-\lambda_{j+1},\,\,j=1,\ldots, n-1, \,\,
 \beta_n=2\lambda_n.
 \end{align*}

In $\mathfrak {osp}_{2s+1|2n}(\mathbb C)_{\bar 0}=\mathfrak{so}_{2s+1}(\mathbb
C)\oplus\mathfrak {sp}_{2n}(\mathbb C)$ we fix the following Cartan subalgebra 
$$
\mathfrak t(\mathfrak {osp}_{2s+1|2n}(\mathbb C)_{\bar 0})=\mathfrak
t_0(\mathfrak {so}_{2s+1}(\mathbb C))\oplus\mathfrak t_1(\mathfrak
{sp}_{2n}(\mathbb C)),
$$
where
 \begin{align*}
\mathfrak t_0(\mathfrak {so}_{2s+1}(\mathbb C))=&
\{\operatorname{diag}(\mu_1,\dots,\mu_{s},0, -\mu_1,\dots,-\mu_{s})\},\\
\mathfrak t_1(\mathfrak {sp}_{2n}(\mathbb C))=&
\{\operatorname{diag}(\lambda_1,\dots,\lambda_n,-\lambda_1,
\dots,-\lambda_n)\}.
 \end{align*}
 We fix the following system of positive roots:
$$
\Delta^+(\mathfrak {osp}_{2s+1|2n}(\mathbb C)_{\bar
	0})=\Delta^+_1(\mathfrak {so}_{2s+1}(\mathbb C)) \cup
\Delta^+_2(\mathfrak {sp}_{2n}(\mathbb C)),
$$
where
\begin{align*}
\Delta^+_1(\mathfrak {so}_{2s+1}(\mathbb C)_{\bar
	0})&=\{\mu_i-\mu_j,\,\mu_i+\mu_j, \,\,i<j, \,\, \mu_i,\},
\\
\Delta^+_2(\mathfrak {sp}_{2n}(\mathbb C))&=\{\lambda_p-\lambda_q,\,\,
\,\,p<q, \,\, \lambda_p+\lambda_q,\,\, \,\,p\leq q \},
\end{align*}
and the following system of simple roots 
$$
\Phi(\mathfrak {osp}_{2s+1|2n}(\mathbb
C)_{\bar 0})=\Phi_1(\mathfrak {so}_{2s+1}(\mathbb C))\cup
\Phi_2(\mathfrak {sp}_{2n}(\mathbb C)),
$$
 where
\begin{align*}
\Phi_1(\mathfrak {osp}_{2s+1}(\mathbb C))&= \{\alpha_1,...,
\alpha_{s}\}, \,\,\, \alpha_i=\mu_i-\mu_{i+1},
\,\,i=1,\ldots, s-1, \,\,\alpha_{s}=\mu_{s},\\
\Phi_2(\mathfrak {sp}_{2n}(\mathbb C))&= \{\beta_1,..., \beta _{n}\},
\,\,\, \beta_j=\lambda_j-\lambda_{j+1},\,\,j=1,\ldots, n-1, \,\,
\beta_n=2\lambda_n.
\end{align*}

In  $\pi\mathfrak {sp}_{n}(\mathbb C)_{\bar 0}=\mathfrak{gl}_n(\mathbb C)$
 we fix the following Cartan subalgebra  
$$
\mathfrak t(\pi\mathfrak {sp}_{n}(\mathbb C)_{\bar 0})=
\{\operatorname{diag}(\mu_1,\dots,\mu_n)\},
$$
 the following system of positive roots:
$$
\Delta^+(\pi\mathfrak {sp}_{n}(\mathbb C)_{\bar 0})=\{\mu_i-\mu_j, \,\,i<j\}
$$
and the following system of simple roots
$$
\begin{array}{c}
\Phi(\pi\mathfrak {sp}_{n}(\mathbb C)_{\bar 0})= \{\alpha_1,..., \alpha_{n-1}\}, \,\,\,
\alpha_i=\mu_i-\mu_{i+1}.
\end{array}
$$

Let $\mathfrak {g} = \mathfrak {osp}_{m|2n}(\mathbb
C)$ or $\pi\mathfrak {sp}_{n}(\mathbb C)$.   Denote by $\mathfrak t(\mathfrak {g}_{\bar 0})^*(\mathbb R)$
the real subspace in $\mathfrak t(\mathfrak {g}_{\bar 0})^*$ 
spanned by ($\mu_j$, $\lambda_i$) or by ($\mu_j$). Consider in $\mathfrak t(\mathfrak {g}_{\bar 0})^*(\mathbb R)$ the scalar product $( \,,\, )$ such that the vectors  $(\mu_j,\lambda_i)$ or $(\mu_j)$ form an orthonormal basis. An element $\gamma\in \mathfrak t(\mathfrak {g}_{\bar 0})^*(\mathbb R)$ is called {\it dominant} if $(\gamma, \alpha)\ge 0$ for all positive roots $\alpha \in \Delta^+(\mathfrak {g}_{\bar 0})$.

Let $G=G^e_{m}$, where $m=2s$ or $2s+1$, or $G^o$. In other words, $G=\OSp_{m|2n}(\mathbb
C)_{\red}$, where $m=2s$ or $2s+1$, or $\Pi\!\Sp_n(\mathbb C)_{\red}$. We also put $P= P^e_{m}$, where $m=2s$ or $2s+1$, or $P^o$, and $R= R^e_{m}$, where $m=2s$ or $2s+1$, or $R^o$.  Let $\mathbf E_{\varphi}\to M= G/P$ be a homogeneous holomorphic vector bundle corresponding to a representation $\varphi$ of $P$ in the fiber $E=(\mathbf E_{\varphi})_{P}$ of $\mathbf E_{\varphi}$ at the point $P$. (An introduction to the theory of homogeneous holomorphic vector bundles can be found in \cite{ADima}.) We denote by $\mathcal E_{\varphi}$ the shaef of holomorphic sections of $\mathbf E_{\varphi}$.

\medskip

\t [Borel-Weil-Bott]. \label{teor borel} {\sl Assume that the representation	$\varphi: P\to \GL(E)$ is completely reducible and $\lambda_1,..., \lambda_s$
	are highest weights of $\varphi|R$. Then the $G$-module $H^0(G/P,\mathcal E_{\varphi})$ is isomorphic to the sum of irreducible $G$-modules with highest weights $\lambda_{i_1},..., \lambda_{i_t}$, where 
	$\lambda_{i_a}$ are dominant highest weights.

}

\medskip

\section{Vector fields on super-Grassmannians}

Let us repeat briefly main definitions from the theory of homogeneous supermanifolds. An action of a Lie supergroup $\mathcal G$ on a supermanifold $\mathcal M$ is called {\it transitive}
if the underlying action of the Lie group $\mathcal G_0$ is transitive on the underlying space $\mathcal M_0$ of $\mathcal M$ and the corresponding  action of the Lie superalgebra is also transitive, see \cite{Viholom} for details.  A supermanifold $\mathcal M$ is called {\it homogeneous} if it possesses a transitive action of a certain Lie supergroup.

 From Formulas (\ref{eq action of Q}) we can deduced that the flag supermanifolds $\mathbf {F}_{k|l}$, $\mathbf {F}^e_{k|l}$ and $\mathbf {F}^o_{k|l}$ are  homogeneous with respect to the action (\ref{eq action of Q}) of the Lie supergroups $\mathcal G = \GL_{m|n}(\mathbb C)$, $\operatorname{OSp}_{m|2n}(\mathbb C)$ and $\Pi\!\Sp_{n}(\mathbb C)$, respectively.  Indeed, the underlying action of $\mathcal G_0$ is transitive, because it is just the standard action of $\mathcal G_0$ on a (isotropic) flag manifold. The corresponding action of the Lie superalgebra is also transitive. For instance this can be explicitly verified using our local coordinates, see also \cite{Viholom}. Further, the actions of Lie supergroups $\operatorname{OSp}_{m|2n}(\mathbb C)$ and $\Pi\!\Sp_{n}(\mathbb C)$ induce the following Lie superalgebra actions
$$
\mu^{e}: \mathfrak {osp}_{m|2n}(\mathbb C)\to
\mathfrak{v}(\mathbf {F}_{k|l}^e)\quad \text{and} \quad \mu^{o}: \pi\mathfrak{sp}_{n}(\mathbb C) \to
\mathfrak{v}(\mathbf {F}_{k|l}^o), 
$$
respectively. We will prove that with some exceptional cases $\mu^{e}$ and $\mu^{o}$ are isomorphisms.

The Lie superalgebras of holomorphic vector fields on isotropic  super-Grass\-mannians corresponding to $\operatorname{OSp}_{m|2n}(\mathbb C)$ and $\Pi\!\Sp_{n}(\mathbb C)$ were calculated in  \cite{oniosp,onipisp}.

\medskip
\t\label{teor super=Grassmannians} {\sl Let $r=1$.

{\bf 1.}	Assume that $m=2k_1$ and $n=2l_1$, i.e the super-Grassmannian\linebreak $\mathbf {F}_{m,k_1|n,l_1}^e$ is of maximal type. If $k_1\geq 1$ and $l_1\geq 1$, then the homomorphism 
$$
\mu^{e}: \mathfrak {osp}_{m|n}(\mathbb C)\to
\mathfrak{v}(\mathbf {F}_{m,k_1|n,l_1}^e)
$$
is in fact an isomorphism. Moreover, if $k_1\geq
2$ and $l_1\geq 1$ the super-Grassmannian $\mathbf {F}_{m-1,k_1-1|n,l_1}^e$  is isomorphic to a connected component of the super-Grassmannian $\mathbf {F}_{m,k_1|n,l_1}^e$.

{\bf 2.} 
Assume that $n=k_1+l_1$, i.e the super-Grassmannian $\mathbf {F}_{n,k_1|n,l_1}^o$ is of maximal type. If $k_1\geq 3$, $l_1\geq 2$, then the homomorphism
$$
\mu^{o}: \pi\mathfrak{sp}_{n}(\mathbb C) \to
\mathfrak{v}(\mathbf {F}_{n,k_1|n,l_1}^o). 
$$
is in fact an isomorphism.

}

\medskip

The Lie superalgebra of holomorphic vector fields on the  super-Grass\-mannian $\mathbf {F}_{m,k_1|n,l_1}$ corresponding to $\mathfrak{gl}_{m|n}(\mathbb C)$ was computed in \cite{buneg,onigl,onich,serov}, see also \cite{ViGL in JA}. The case of flag supermanifolds corresponding to $\mathfrak{gl}_{m|n}(\mathbb C)$ was studied in  \cite{ViGL in JA}. The case of super-Grassmannian  corresponding to $\mathfrak{q}_{n}(\mathbb C)$ can be found in \cite{oniq}, and the case of flag supermanifolds of this type can be found in \cite{ViPi-sym}.

\section{Vector fields on isotropic flag supermanifolds}

Assume that $r>1$. From now on we consider isotropic supermanifolds of maximal type.  As we have seen in Section $4$ the isotropic flag supermanifold $ \mathcal M^e: =\mathbf {F}_{k|l}^e$ (or $\mathcal M^o:=\mathbf {F}_{k|l}^o$) is a superbundle. We denote by $\mathcal B^e$ and $\mathcal S^e$ (or by $\mathcal B^o$ and $\mathcal S^o$) its base space and its fiber, respectively. In other words we put 
$$
\mathcal M^e=\mathbf {F}_{k|l}^e, \quad \mathcal B^e= 
\mathbf {F}_{k_0,k_1|l_0,l_1}^e\quad \text{and} \quad  \mathcal S^e=\mathbf{F}_{k'|l'}
$$
for $\mathfrak {osp}_{m|2n}(\mathbb C)$-flag supermanifolds, and
$$
\mathcal M^o=\mathbf {F}_{k|l}^o, \quad \mathcal B^o= 
\mathbf {F}_{k_0,k_1|l_0,l_1}^o\quad \text{and} \quad  \mathcal S^o=\mathbf{F}_{k'|l'}
$$
for $\pi\mathfrak{sp}_{n}(\mathbb C)$-flag supermanifolds. Note that in both cases the fiber is a usual $\mathfrak{gl}_{k_1|l_1}(\mathbb C)$-flag supermanifold.

\subsection{Plan of the proof} Consider the case of $\mathfrak {osp}_{2m|2n}(\mathbb C)$-flag supermanifolds.  Assume that we do not have non-constant holomorphic functions on the fiber $\mathcal S^e$, i.e.  $\mathcal O_{\mathcal S^e}(\mathcal S^e_0)=\mathbb C$. Then by Proposition \ref{prop Bash}, the projection
$\mathcal M^e\to \mathcal B^e$ determines the homomorphism of Lie superalgebras $\mathcal{P}^e:\mathfrak v(\mathcal M^e)\to\mathfrak
v(\mathcal B^e).$
This projection is $\operatorname{OSp}_{2m|2n}(\mathbb C)$-equivariant. Hence for the Lie superalgebra homomorphisms 
$$
\mu^e: \mathfrak {osp}_{2m|2n}(\mathbb C)\to\mathfrak
v(\mathcal M^e)\quad \text{and} \quad \mu^e_{\mathcal B}: \mathfrak {osp}_{2m|2n}(\mathbb C)\to\mathfrak v(\mathcal B^e)
$$
 we have $
\mu^e_{\mathcal B} = \mathcal{P}^e\circ\mu^e.$ Assuming conditions of Theorem \ref{teor super=Grassmannians}, the homomorphisms $\mu^e_{\mathcal B}$ and hence the homomorphism
$\mathcal{P}^e$ is surjective.  In case if $\mathcal{P}^e$ is also injective we have
\begin{equation}\label{eq connection mu and mu_B}
\mu^e = (\mathcal{P}^e)^{-1}\circ\mu^e_{\mathcal B}
\end{equation}
is surjective and therefore  $\mathfrak v(\mathcal M^e)\simeq \mathfrak {osp}_{2m|2n}(\mathbb C).$ Hence our goal is to prove that $\mathcal{P}^e$ is injective. For $\pi\mathfrak{sp}_{n}(\mathbb C)$-flag supermanifolds the idea is similar.

In the case of $\mathfrak {osp}_{2m-1|2n}(\mathbb C)$-flag supermanifolds a similar argument does not work since $\mu^e_{\mathcal B}: \mathfrak {osp}_{2m-1|2n}(\mathbb C)\to\mathfrak v(\mathcal B^e)$ is not surjective, see Theorem \ref{teor super=Grassmannians}. This case we will consider in a separate paper.

\subsection{Vector bundles $\mathbf W^e_0$ and $\mathbf W^o_0$} 

 In  Section $4$ we defined the locally free sheaf $\tilde{\mathcal W}_0$ for any superbundle $\mathcal M$. Denote by $\tilde{\mathcal W}_0^e$ and $\tilde{\mathcal W}_0^o$  the locally free sheaves corresponding to the superbundles $\mathbf {F}_{k|l}^e$ and $\mathbf {F}_{k|l}^o$, respectively. We denote also by $\mathbf W_0^e$ and $\mathbf W_0^o$ the corresponding to $\tilde{\mathcal W}_0^e$ and $\tilde{\mathcal W}_0^o$ vector bundles  over $\mathcal B^e_0 = G_m^e/P_m^e$ and $\mathcal B^o_0= G^o/P^o$, respectively.  Our goal now is to compute the vector space of global sections of $\mathbf W_0^e$ and $\mathbf W_0^o$  using  Theorem \ref{teor borel}. Note that the vector bundles $\mathbf W_0^e$ and $\mathbf W_0^o$ are homogeneous because the sheaves  $\tilde{\mathcal W}_0^e$ and $\tilde{\mathcal W}_0^o$ possess the natural actions of the Lie groups $\operatorname{OSp}_{m|2n}(\mathbb C)_{\red}$ and $\Pi\!\Sp_{n}(\mathbb C)_{\red}$, respectively. 
Let us compute the corresponding to $\mathbf W_0^e$ and $\mathbf W_0^o$ representations of $P_m^e$ and $P^o$.

\subsection{\bf Representations of $P_{2m}^e$  and $P^o$}

Consider the local chart on the
super-Grassmannian $\mathcal B^e$ of maximal type corresponding to 
\begin{equation}\label{eq I_1 of o in  B_0}
I^e_{1\bar 0} = \{m-k_1+1,\ldots,m\}\quad \text{and}\quad I^e_{1\bar 1} = \{n-l_1+1,\ldots,n\},
\end{equation}
where $m=2k_1$ and $n=2l_1$, 
and the local chart on the
super-Grassmannian $\mathcal B^o$  of maximal type corresponding to 
\begin{equation}\label{eq I_1 of o in  B_0o}
I^o_{1\bar 0} = \{n-k_1+1,\ldots,m\}\quad \text{and}\quad I^o_{1\bar 1} = \{1,\ldots,l_1\},
\end{equation}
where $k_1+l_1=n$. We put $I_1^e=(I^e_{1\bar 0},I^e_{1\bar 1})$ and $I_1^o=(I^o_{1\bar 0},I^o_{1\bar 1})$. The coordinate matrices $ Z_{I_1^e}$ and $ Z_{I_1^o}$ in this cases have the following form, respectively. 
\begin{equation}\label{eq loc chart on B}
 {\bf 1.}\,\,
 Z_{I_1^e} =\left(
 \begin{array}{cc}
 X_1&\Xi_1\\
 E_{k_1}&0\\
 -\Xi_1^t&Y_1\\0&E_{l_1}\end{array} \right),\quad
 {\bf 2.}\,\,
 Z_{I_1^o} =\left(
 \begin{array}{cc}
 X_1&\Xi_1\\
 E_{k_1}&0\\
 0&E_{l_1}\\
 \H_1&-X^t_1\end{array} \right).
\end{equation}
 {\bf 1.} Here $X_1$, $Y_1$ are matrices of size $k_1\times k_1$ and $l_1\times l_1$, respectively, that contain even coordinates, and $\Xi_1$ is a matrix of size $k_1\times l_1$, that contains odd coordinates. Moreover we have $X_1^t=-X_1$ and $Y_1^t=Y_1$. 
{\bf 2.} $X_1$ is a matrix of size $l_1\times k_1$ that contains even coordinates, and $\Xi_1$, $\H_1$  are matrices of size $l_1\times
l_1$ and $k_1\times k_1$, respectively, that contain odd coordinates. Moreover we have $\Xi_1^t=\Xi_1$, $\H_1^t=-\H_1$.

Denote by $x^e$ the point in $\mathcal B^e_0$ defined by the following equations
$$
X_1=0 ,\quad Y_1= 0, \quad \Xi_1=0,
$$
and by $x^o$ the point in $\mathcal B^o_0$ defined by the following equations
$$
X_1=0 ,\quad \Xi_1= 0, \quad \H_1=0.
$$
It is easy to see that the Lie groups $P^e_{2m}$ and $P^o$, see (\ref{eq parabolic OSP})  and (\ref{P_pisp}), are stabilizers of $x^e$ and $x^o$, respectively. Recall that we denoted by  $R^e_{2m}$ and $R^o$ the reductive parts of $P^e_{2m}$ and $P^o$, respectively.

 Let us compute the representation $\psi^e$ and $\psi^o$ of $P^e_{2m}$ and $P^o$ in the fibers $(\mathbf W^e_0)_{P_{2m}^e}$ and $(\mathbf W^o_0)_{P^o}$, respectively. We identify $(\mathbf W^e_0)_{P_{2m}^e}$ (or $(\mathbf W^o_0)_{P^o}$) with the Lie superalgebra of holomorphic vector fields $\mathfrak v(\mathcal S^e)$ (or $\mathfrak v(\mathcal S^o)$), see Proposition \ref{prop tilde W_0 is localy free}. 
Let us choose an atlas on $\mathcal M^e$ and $\mathcal M^o$  defined by $I^e_{1}$ and $I^o_{1}$, as above, and by certain $I^e_{s},\; s = 2,\ldots,r$ and $I^o_{s},\; s = 2,\ldots,r$. In notations (\ref{eq loc chart on B})  the Lie groups $P^e_{2m}$ and $P^o$ act at $x^e$ and $x^o$ in the chart on the super-Grassmannians $\mathcal B^e$ and $\mathcal B^o$ defined by $Z_{I^e_1}$ and $Z_{I^o_1}$ in the following way: 
$$
\begin{array}{l}
{\bf 1.}\quad\left(
\begin{array}{cccc} A_{1}&0&0
&0\\C_1&(A_{1}^T)^{-1}&0&0\\0&0&A_{2}&0\\0&0&C_2&(A_{2}^T)^{-1}\end{array}
\right) \left(
\begin{array}{cc}
0&0\\
E_{k_1}&0\\
0&0\\0&E_{l_1}\end{array} \right) = \left(
\begin{array}{cc} 0&0\\  (A_1^T)^{-1} & 0\\
0 & 0\\ 0 &  (A_2^T)^{-1}
\end{array} \right),\\
{\bf 2.}\quad\left(
\begin{array}{cccc}
A_{1}&0&0
&0\\C_1&B_{1}&0&0\\0&0&(A_{1}^{T})^{-1}&D_1\\
0&0&0&(B_{1}^T)^{-1}
\end{array}
\right) \left(
\begin{array}{cc}
0&0\\
E_{k_1}&0\\
0&E_{l_1}\\
0&0\end{array} \right)
= \left(
\begin{array}{cc} 0&0\\  B_1 & 0\\
0& (A_1^{-1})^T\\
0& 0
\end{array} \right).
\end{array}
$$
(Note that a chart on $\mathcal B^e$ is defined by $Z_{I^e_1}$, and a chart on the whole flag supermanifold $\mathcal M^e$ is defined by $Z_{I^e_s}$, where $s= 1,\ldots, r$. The same holds for $\mathcal M^o$.) Further, for $Z_{I^e_2}$ and $Z_{I^o_2}$, we have
\begin{equation}
\label{eq action over o}
\begin{split}
&{\bf 1.} \quad\left(
\begin{array}{cc} (A_1^T)^{-1} & 0\\ 0 &  (A_2^T)^{-1}\end{array}
\right) \left(
\begin{array}{cc}  X_2 & \Xi_2 \\ \H_2 & Y_2\end{array} \right) =
\left(
\begin{array}{cc}  (A_1^T)^{-1}X_2  &
(A_1^T)^{-1}\Xi_2 \\
(A_2^T)^{-1} \H_2 &
(A_2^T)^{-1}Y_2
\end{array} \right),\\
&\\&
{\bf 2.} \quad\left(
\begin{array}{cc} B_1 & 0\\0 &  (A_1^T)^{-1}\end{array}
\right) \left(
\begin{array}{cc}  X_2 & \Xi_2 \\ \H_2 & Y_2\end{array} \right) =
\left(
\begin{array}{cc}  B_1X_2  &
B_1\Xi_2 \\
(A_1^T)^{-1} \H_2 & (A_1^T)^{-1}Y_2
\end{array} \right).
\end{split}
\end{equation}
Note that the local coordinates of $Z_{I^e_s},\; s \ge 2$, (or $Z_{I^o_s},\; s \ge 2$) can be interpreted as local coordinates on the fiber $\mathcal S^e$ (or $\mathcal S^o$) of the superbundle $\mathcal M^e$ (or $\mathcal M^o$). Hence to obtain the actions of $P^e_{2m}$ and $P^o$ we use (\ref{eq action over o}) and modify
$Z_{I^e_s}$ and $Z_{I^o_s}$, $s \ge 3$, accordingly. We see that the nilradicals of $P^e_{2m}$ and $P^o$ act trivially on
$\mathcal S^e$ and $\mathcal S^o$. Further, the action of the reductive parts $R^e_{2m}$ and $R^o$ coincide with the restriction 
\begin{description}
\item[1.] of $R^e_{2m}$ on the subgroup $\GL_{k_1}(\mathbb C)\times \GL_{l_1}(\mathbb C)$ that acts on $\mathcal S^e$ as the composition of the standard action of
$\GL_{k_1|l_1}(\mathbb C)_{\red}$, see (\ref{eq action of Q}), and 
$$
\GL_{k_1}(\mathbb C)\times \GL_{l_1}(\mathbb C)\to \GL_{k_1}(\mathbb C)\times
\GL_{l_1}(\mathbb C),\,\, (A_1,A_2)\mapsto ((A_1^T)^{-1},(A_2^T)^{-1});
$$

\item[2.] of $R^o$ on the subgroup $\GL_{k_1}(\mathbb C)\times \GL_{l_1}(\mathbb C)$ that acts on $\mathcal S^o$ as a composition of the standard action of
$\GL_{k_1|l_1}(\mathbb C)_{\red}$, see (\ref{eq action of Q}), and 
$$
\GL_{k_1}(\mathbb C)\times \GL_{l_1}(\mathbb C)\to \GL_{k_1}(\mathbb C)\times
\GL_{l_1}(\mathbb C),\,\, (B_1,A_1)\mapsto (B_1,(A_1^T)^{-1}).
$$
\end{description}

Assume that 
$$
\mathfrak v(\mathcal S^e)\simeq \mathfrak v(\mathcal S^o)\simeq \mathfrak{gl}_{k_1|l_1}(\mathbb C)/ \langle
E_{k_1+l_1}\rangle = \left\lbrace\left(
\begin{array}{cc}
Z_1 & T_1  \\
T_2 & Z_2  \\
\end{array}
\right)+<E_{k_1+l_1}>\right\rbrace,
$$
where $Z_1\in \mathfrak{gl}_{k_1}(\mathbb C)$ and $Z_2\in \mathfrak{gl}_{l_1}(\mathbb C)$. Then 
\begin{description}
 \item[1.] the representation $\psi^e$ of $R^e_{2m}$ on  $\mathfrak v(\mathcal S^e)$ is determined by 
$$
\begin{array}{c}
\left(
\begin{array}{cc}
(A_1^T)^{-1} & 0 \\
0 & (A_2^T)^{-1}
\end{array}
\right) \left(                              \left(
\begin{array}{cc}
Z_1 & T_1  \\
T_2 & Z_2  \\
\end{array}
\right)+<E_{k_1+l_1}>\right)
\left(
\begin{array}{cc}
A_1^T & 0 \\
0 & A_2^T
\end{array}
\right)=\\
\left(
\begin{array}{cc}
(A_1^T)^{-1}Z_1A_1^T & (A_1^T)^{-1}T_1A_2^T \\
(A_2^T)^{-1}T_1A_1^T & (A_2^T)^{-1}Z_2A_2^T
\end{array}
\right)+<E_{k_1+l_1}>;
\end{array}
$$
where $A_1\in \GL_{k_1}(\mathbb C)$, $A_2\in \GL_{l_1}(\mathbb C)$;

\item[2.] the representation $\psi^o$ of $R^o$ on  $\mathfrak v(\mathcal S^o)$ is determined by 
$$
\begin{array}{c}
\left(
\begin{array}{cc}
B_1 & 0 \\
0 & (A_1^T)^{-1}
\end{array}
\right) \left(                              \left(
\begin{array}{cc}
Z_1 & T_1  \\
T_2 & Z_2  \\
\end{array}
\right)+<E_{k_1+l_1}>\right)
\left(
\begin{array}{cc}
B_1^{-1} & 0 \\
0 & A_1^T
\end{array}
\right)=\\
\left(
\begin{array}{cc}
B_1Z_1B_1^{-1} & B_1T_1A_1^T \\
(A_1^T)^{-1}T_1B_1^{-1} & (A_1^T)^{-1}Z_2A_1^T
\end{array}
\right)+<E_{k_1+l_1}>,
\end{array}
$$
where $B_1\in \GL_{k_1}(\mathbb C)$, $A_1\in \GL_{l_1}(\mathbb C)$.
\end{description}

Denote by $\rho_{1}$ and $\rho_{2}$ the standard representations of $\GL_{k_1}(\mathbb C)$ and $\GL_{l_1}(\mathbb C)$ on $\mathbb {C}^{k_1}$ and $\mathbb {C}^{l_1}$, respectively, and by $\Ad_{1}$ and $\Ad_{2}$ the adjoint representations of $\GL_{k_1}(\mathbb C)$ and $\GL_{l_1}(\mathbb C)$ on $\mathfrak {sl}_{k_1}(\mathbb C)$ and $\mathfrak {sl}_{l_1}(\mathbb C)$, respectively. We also denote by
$1$ the one dimensional trivial representation of $\GL_{k_1}(\mathbb C)\times \GL_{l_1}(\mathbb C)$.
We have proved the following lemma.

\medskip

\lem\label{lem representation of Gl_k_1 in v(S)} 
 {\sl The representations $\psi^e$ and $\psi^o$ of $P^e_{2m}$ and $P^o$ in the fibers $(\mathbf W^e_0)_{P_{2m}^e}$ and $(\mathbf W^o_0)_{P^o}$, respectively, are completely reducible.  If
$\mathfrak v(\mathcal S)\simeq \mathfrak{gl}_{k_1|l_1}(\mathbb C) / \langle
E_{k_1+l_1}\rangle $, 
then
\begin{description}
	\item[1.] \begin{equation}
	\begin{split}
	\label{predstavl_v_W_0_osp}  
	\psi^e|R^e_{2m} =\left\{
	\begin{array}{l}
	\Ad_{1} + \Ad_{2} + \rho^*_{1}\otimes\rho_{2} +
	\rho_{1}\otimes\rho^*_{2} + 1\; \text{for}\; k_1,l_1 > 0,\\
	\Ad_{1}\; \text{for}\; k_1>0,\; l_1 = 0,\\
	\Ad_{2}\; \text{for}\; k_1 = 0,\; l_1 > 0.\end{array}\right.
	\end{split}
	\end{equation}
	\item[2.] \begin{equation}\begin{split}
	\label{predstavl_v_W_0 pisp}
	\psi^o|R^o
	=\left\{
	\begin{array}{l}
	\Ad_{1} + \Ad_{2} + \rho_{1}\otimes\rho_{2} +
	\rho^*_{1}\otimes\rho^*_{2} + 1\; \text{for}\; k_1,l_1 > 0,\\
	\Ad_{1}\; \text{for}\; k_1>0,\; l_1 = 0,\\
	\Ad_{2}\; \text{for}\; k_1 = 0,\; l_1 >
	0.\end{array}\right.
\end{split}
	\end{equation}
\end{description}

}

\medskip

Further we will use the charts $\mathcal U^e$ and $\mathcal U^o$ on $\mathbf {F}^e_{k|l}$ and $\mathbf {F}^o_{k|l}$
defined by 
$I^e_s= I^e_{s\bar 0}\cup I^e_{s\bar 1}$ and $I^o_s= I^o_{s\bar 0}\cup I^o_{s\bar 1}$, where $I^e_{1\bar i}$ and $I^o_{1\bar i}$ are as above, and  
$$
I^e_{s\bar 0}=I^o_{s\bar 0}=\{k_{s-1}-k_s+1,\ldots,k_{s-1}\},\quad I^e_{s\bar 1}=I^o_{s\bar 1}= \{l_{s-1}-l_s+1,\ldots,l_{s-1}\}
$$
for $s\ge 2$.
The coordinate matrices of these charts have the following form
$$
Z_{I^e_s} = Z_{I^o_s} = \left(
\begin{array}{cc} 
X_s&\Xi_s\\ 
E_{k_s}&0\\
\H_s&Y_s\\
0&E_{l_s}
\end{array}\right),\quad s=2,\ldots, r,
$$
where again the local coordinates are $X_s=(x^s_{ij})$, $Y_s=(y^s_{ij})$,  $\Xi_s=(\xi^s_{ij})$ and $\H_s=(\eta^s_{ij}).$
We denote by $\mathcal U^e_{\mathcal B^e}$ and by $\mathcal U^o_{\mathcal B^o}$ the corresponding charts on $\mathcal B^e$ and $\mathcal B^o$, respectively. In other words, $\mathcal U^e_{\mathcal B^e}$ and $\mathcal U^o_{\mathcal B^o}$ are given by  (\ref{eq loc chart on B}).

The proofs of the following two lemmas are similar to the proof of Lemma $2$ and Lemma $3$ in \cite{ViGL in JA}.

\medskip

\lem \label{lem some fundamental vector fields} {\sl The following vector fields ${\bf 1.}\,\,\,\frac{\partial}{\partial \xi^1_{ij}}$, ${\bf 2.}
	\,\,\,\frac{\partial}{\partial \xi^1_{ij}},\,\,
	\frac{\partial}{\partial \eta^1_{ij}}$ 	 are fundamental. That is they are induced by the natural actions of Lie supergroups $\operatorname{OSp}_{2m|2n}(\mathbb C)$ and $\Pi\!\Sp_{n}(\mathbb C)$  on $\mathbf {F}^e_{k|l}$ and $\mathbf {F}^o_{k|l}$, respectively.
}

\medskip

\lem \label{lem assume ker Pi ne 0} {\sl Assume that $\Ker \mathcal{P}\ne \{0\}$. Then $\dim \mathcal{W}^e_{(0)}(\mathcal B^e_0)> \dim  \mathcal{W}^e_{(1)}(\mathcal B^e_0)$ and $\dim \mathcal{W}^o_{(0)}(\mathcal B^o_0) > \dim  \mathcal{W}^o_{(1)}(\mathcal B^o_0)$.	}

\medskip

Here $\mathcal{W}^e$ and $\mathcal{W}^o$ are locally free sheaves as in Section $4$ corresponding to superbundles $\mathbf {F}^e_{k|l}$ and $\mathbf {F}^o_{k|l}$, respectively.

Now we need the following lemma. 

\medskip

\lem \label{lem sections of W_0} {\sl Assume  that $\mathcal O_{\mathcal S^e}(\mathcal S^e_0) \simeq \mathbb C$,  $\mathcal O_{\mathcal S^o}(\mathcal S^o_0) \simeq \mathbb C$,
$
\mathfrak v(\mathcal S^e)\simeq \mathfrak{pgl}_{k_1|l_1}(\mathbb C)$ and $\mathfrak v(\mathcal S^o)\simeq \mathfrak{pgl}_{k_1|l_1}(\mathbb C).
$
 Then 
 $$
\tilde{\mathcal W}_0(\mathcal B^e_0)\simeq \left\{
 \begin{array}{ll}
 \mathbb C, & k_1>2,\,l_1\geq 1; \\
 \mathbb C\oplus \mathfrak{r}_1, & k_1=2,\,l_1\geq 1; \\
 \mathbb C\oplus \mathfrak{r}_2, & k_1=1,\,l_1> 1, \\
 \end{array}
 \right.
 $$
 where $\mathfrak{r}_1$ is the $\mathfrak{so}_{4}(\mathbb C)\oplus
 \mathfrak{sp}_{n}(\mathbb C)$-module with the highest weight  $\mu_{1}-\mu_{2}$,
 $\mathfrak{r}_2$ is the  $\mathfrak {sp}_{2l_1}(\mathbb C)_{\bar
 	0}$-module with the highest weight $\lambda_{1}$ and $\mathbb C$ is the trivial
 $\mathfrak{so}_{m}(\mathbb C)\oplus \mathfrak{sp}_{n}(\mathbb C)$-module that corresponds to the highest weight $0$.
 
Further we have
 $$
\tilde{\mathcal W}_0(\mathcal B^o_0)\simeq 
 \left\{
 \begin{array}{ll}
 \mathbb C, & k_1>1,\,l_1> 1; \\
 \mathbb C\oplus \mathfrak{r}_1, & k_1=1,\,l_1> 1; \\
 \mathbb C\oplus \mathfrak{r}_2, & k_1>1,\,l_1= 1, \\
 \end{array}
 \right.
 $$
 where $\mathfrak{r}_1$ is the $\mathfrak{sl}_{n}(\mathbb C)$-module with the highest weight $-\mu_{l_1}- \mu_{l_1+1}$, $\mathfrak{r}_2$ is the
 $\mathfrak{sl}_{n}(\mathbb C)$-module with the highest weight $\mu_{1}+\mu_{2}$,
 $\mathbb C$ is the trivial $\mathfrak{sl}_{n}(\mathbb C)$-module that corresponds to the highest weight $0$.

}
\medskip

\noindent{\it Proof.} {\bf 1.} We compute the vector space of global sections of $\mathbf W^e_0$ using Theorem \ref{teor borel}. The representation $\psi^e$ of the Lie group $P^e_{2m}$ was computed in Lemma  \ref{lem representation of Gl_k_1 in v(S)}.
It follows that the highest weights of $\psi^e$  have the following form:
\begin{itemize}
	\item $\mu_{1}-\mu_{k_1}$, $\mu_{1}-\lambda_{l_1}$,
	$\lambda_{1}- \mu_{k_1}$, $\lambda_{1}- \lambda_{l_1}$,
	$0$ for $k_1>1$, $l_1>1$;
	\item $\mu_{1}-\lambda_{l_1}$,
	$\lambda_{1}- \mu_{1}$, $\lambda_{1}- \lambda_{l_1}$,
	$0$ for $k_1=1$, $l_1>1$;
	\item $\mu_{1}-\mu_{k_1}$, $\mu_{1}-\lambda_{1}$,
	$\lambda_{1}- \mu_{k_1}$,
	$0$ for $k_1>1$, $l_1=1$;
	\item $\mu_{1}-\lambda_{1}$,
	$\lambda_{1}- \mu_{1}$,
	$0$ for $k_1=1$, $l_1=1$.
\end{itemize}
(Note that we have  $k_1>0$ and $l_1>0$, since otherwise $m=0$ or $n=0$.) Therefore the dominant weights of  $\psi^e$ are:
\begin{itemize}
	\item $0$, if $k_1>2$, $l_1\geq 1$;
	
	\item $0$, $\mu_{1}-\mu_{k_1}$, if $k_1=2$, $l_1\geq 1$;
	
	\item $0$, $\lambda_{1}-\mu_{1}$, if $k_1=1$, $l_1>1$;
	
	\item $0$, $\lambda_{1}-\mu_{1}$, if $k_1=1$, $l_1=1$.
	This case we will not consider further since if $k_1=1$, $l_1=1$ we have  $\mathcal O_{\mathcal S^e}(\mathcal S^e_0) \ne \mathbb C$, see Theorem \ref{teor constant functions}.
	
\end{itemize}
By the Borel-Weil-Bott Theorem we get the result.

\medskip

{\bf 2.}  Again we use the Borel-Weil-Bott Theorem to compute the vector space of global sections of $\mathbf W^o_0$. The representation $\psi^o$ of the Lie group $P^o$ was computed in Lemma \ref{lem representation of Gl_k_1 in v(S)}.
It follows that the highest weights of $\psi^o$  have the following form:
\begin{itemize}
	\item $\mu_{l_1+1}-\mu_{n}$, $\mu_{l_1+1}+\mu_{1}$,
	$-\mu_{n}- \mu_{l_1}$, $\mu_{1}- \mu_{l_1}$,
	$0$ for $k_1>1$, $l_1>1$;
	
	\item $\mu_{l_1+1}+\mu_{1}$,
	$-\mu_{l_1}- \mu_{l_1+1}$, $\mu_{1}- \mu_{l_1}$,
	$0$ for $k_1=1$, $l_1>1$;
	
	\item $\mu_{2}-\mu_{n}$, $\mu_{2}+\mu_{1}$,
	$-\mu_{n}- \mu_{1}$,
	$0$ for $k_1>1$, $l_1=1$;
	
	\item $\mu_{1}+\mu_{2}$,
	$-\mu_{2}- \mu_{1}$,
	$0$ for $k_1=1$, $l_1=1$.
\end{itemize}
(Note that by definition of $\mathbf F_{k|l}^o$
we have $k_1>0$ and $l_1>0$. Indeed, if for example $k_1=0$, then $l_1=n$ and $\mathbf F_{k|l}^o$ is isomorphic to a point.)
Therefore the dominant weights of  $\psi^o$ are:
\begin{itemize}
	\item $0$, if $k_1>1$, $l_1>1$;
	
	\item $0$, $-\mu_{l_1}- \mu_{l_1+1}$, if $k_1=1$, $l_1>1$;
	
	\item $0$, $\mu_{1}+\mu_{2}$, if $k_1>1$, $l_1=1$;
\end{itemize}
	The case $k_1=1$ and $l_1=1$ we will not consider further since in this case we have  $\mathcal O_{\mathcal S^o}(\mathcal S^o_0) \ne \mathbb C$, see Theorem \ref{teor constant functions}.
By the Borel-Weil-Bott Theorem we get the result.$\Box$

\medskip

\subsection{Main results}

Now we are ready to prove the following two theorems.

\medskip

\t \label{teor osp} {\it Assume that $r>1$, $m=k_1$ and $n=l_1$. If $\mathcal O_{\mathcal S^e}(\mathcal S^e_0) \simeq \mathbb C$, $\mathfrak v(\mathbf
	F_{2m,k_1|2n,l_1}^e)\simeq \mathfrak
	{osp}_{2m|2n}(\mathbb C)$ and $\mathfrak v(\mathcal S^e)\simeq \mathfrak {pgl}_{k_1|l_1}(\mathbb C)$,
	then $\mathfrak v(\mathbf F_{k|l}^e)\simeq
	\mathfrak {osp}_{2m|2n}(\mathbb C)$.}

\medskip 

\noindent{\it Proof.}  Consider the super-stabilizer $\mathcal P^e_{2m}\subset \operatorname{OSp}_{2m|2n}(\mathbb C)$ of $x^e$. It contains all super-matrices of the following form:
\begin{equation}\label{eq P^e contais matrices super case}
\left(
\begin{array}{cccc} 
A_{1}&0&C_{11} & 0\\
C_1&(A_{1}^T)^{-1}&C_{21} &  C_{22}\\
-C^T_{22} & 0 &A_{2}&0\\
 C^T_{21} &  C^T_{11}&C_2&
(A_{2}^T)^{-1}
\end{array}
\right),
\end{equation} 
where the size of the matrices is as in (\ref{eq parabolic OSP}). Denote by $\mathcal L$ the Lie subsupergroup in $\mathcal P^e_{2m}$ defined by the following  submatrix. 
$$
\left(
\begin{array}{cc}
(A_{1}^T)^{-1} & C_{22} \\
C^T_{11} & (A_{2}^T)^{-1}
\end{array}
\right).
$$
We see that $\mathcal L\simeq \GL_{k_1|l_1}(\mathbb C)$. And if we replace $(A_{1}^T)^{-1}$ by $W_1\in \GL_{k_1}(\mathcal C)$ and $(A_{2}^T)^{-1}$ by $W_2\in \GL_{l_1}(\mathcal C)$ we will see that $\mathcal L$ acts $(\mathbf W^e_0)_{x^e} \simeq \mathfrak {pgl}_{k_1|l_1}(\mathbb C)$ in the standard way. In other words, $(\mathbf W^e_0)_{x^e} \simeq \mathfrak {pgl}_{k_1|l_1}(\mathbb C)$  is isomorphic to the adjoint $\mathfrak l$-module, where $\mathfrak l\simeq \mathfrak {gl}_{k_1|l_1}(\mathbb C)$ is the Lie superalgebra of $\mathcal L$.

Now we repeat the argument used in \cite{ViGL in JA}. 
Let $\pi : \mathcal W^e\to \widetilde{\mathcal W}^e_0 = \mathcal W^e/\mathcal W^e_{(1)}$ be the natural map and $\pi_{x^e} : \mathcal W^e\to (\mathbf W^e_0)_{x^e}$ be the composition of $\pi$ and of the evaluation map at the point $x^e$. We have the following commutative diagram: 
$$
\begin{CD}
\mathcal W^e(\mathcal B_0)@>{[X,\,\,\cdot
	\,\,]}>> \mathcal W^e(\mathcal B_0)
\\
@V{\pi_{x^e}}VV @V{\pi_{x^e}}VV\\
(\mathbf W^e_0)_{x^e}@>{[X,\,\,\cdot
	\,\,]}>> (\mathbf W^e_0)_{x^e}
\end{CD},
$$
where $X\in \mathfrak{l}$. (Note that the vector space $\mathcal W^e(\mathcal B_0)$ is an ideal in $\mathfrak v(\mathcal M^e)$ and in particular it is invariant with respect to the action of $\mathcal L$.)
Denote by $V$ the image $\pi_{x^e}(\mathcal W^e(\mathcal B_0))$. From the commutativity of this diagram it follows that 
$
V\subset (\mathbf W^e_0)_{x^e} \simeq \mathfrak
{pgl}_{k_1|l_1}(\mathbb C)
$
 is invariant with respect to the adjoint representation of $\mathfrak {pgl}_{k_1|l_1}(\mathbb C)$. Therefore, $V$ is an ideal in  $\mathfrak {pgl}_{k_1|l_1}(\mathbb C)$. 
 
 Analyzing ideals in  $\mathfrak {pgl}_{k_1|l_1}(\mathbb C)$, we see that $V\subset \Im(\gamma)$, where $\gamma : \tilde{\mathcal W}^e_0(\mathcal B_0) \to (\mathbf W^e_0)_{x^e}$ is the evaluation map,  never coincides with non-trivial ideals. Hence, $V=\{0\}$.
In other words, we proved that all sections of  $\pi(\mathcal W^e(\mathcal B_0))$ are equal to $0$ at the point $x^e$. Since $\mathbf W^e_{0}$ is a homogeneous bundle, we get that sections from $\pi(\mathcal W^e(\mathcal B_0))$ are equal to $0$ at any point. Therefore, we have $\pi(\mathcal W^e(\mathcal B_0))=\{0\}$ and $\mathcal W^e(\mathcal B_0)_{(0)} \simeq \mathcal W^e(\mathcal B_0)_{(1)}.$
 From Lemma  \ref{lem assume ker Pi ne 0}, it follows that $\Ker\mathcal{P}=\{0\}$.$\Box$

\medskip

We have proved the following theorem.

\medskip

\t \label{teor general osp} {\sl Assume that $r>1$, $m=k_1$, $n=l_1$; the conditions (\ref{eq_condition no functions}) hold; $k_1\geq 1$, $l_1\geq 1$ and $\mathfrak{v}(\mathbf{F}_{k'|l'})\simeq \mathfrak {pgl}_{k_1|l_1}(\mathbb
	C)$. Then
	$\mathfrak
	v(\mathbf F_{k|l}^e)\simeq 	\mathfrak {osp}_{2m|2n}(\mathbb C)$.$\Box$}

\medskip

\medskip

\t \label{teor pisp} {\it Assume that $r>1$ and $n=k_1+l_1$. If $\mathcal O_{\mathcal S^o}(\mathcal S^o_0) \simeq \mathbb C$, $\mathfrak v(\mathbf
	F_{n,k_1|n,l_1}^o)\simeq \pi\mathfrak
	{sp}_{n}(\mathbb C)$ and $\mathfrak v(\mathcal S^o)\simeq \mathfrak
	{pgl}_{k_1|l_1}(\mathbb C)$, then
	$\mathfrak v(\mathbf F_{k|l}^o)\simeq \pi\mathfrak
	{sp}_{n}(\mathbb C)$.

	}

\medskip 

\noindent{\it Proof.}  Consider the super-stabilizer $\mathcal P^o\subset \Pi\!\Sp_{n}(\mathbb C)$ of $x^o$. It contains all super-matrices of the following form:
\begin{equation}\label{eq P^o contais matrices super case}
\left(
\begin{array}{cccc} 
A_{1}&0&0 &  C_{12}\\
C_1&B_{1}&C_{12}^T & C_{22}\\
 D_{11} &  D_{12} &(A_{1}^{T})^{-1}&-C_1\\
- D^T_{12} & 0&0&(B_{1}^T)^{-1}
\end{array}
\right),
\end{equation} 
where the size of the matrices is as in (\ref{P_pisp}). Denote by $\mathcal L$ the Lie subsupergroup in $\mathcal P$ defined by the following coordinate matrix. 
$$
\left(
\begin{array}{cc}
B_1 & C_{22} \\
D_{12} &(B_{1}^T)^{-1}
\end{array}
\right).
$$
We see that $\mathcal L\simeq \GL_{k_1|l_1}(\mathbb C)$. And if we replace $(B_{1}^T)^{-1}$ by $W_1\in \GL_{l_1}(\mathcal C)$  we will see that $\mathcal L$ acts $(\mathbf W^o_0)_{x^o} \simeq \mathfrak {pgl}_{k_1|l_1}(\mathbb C)$ in the standard way. In other words, $(\mathbf W^o_0)_{x^o} \simeq \mathfrak {pgl}_{k_1|l_1}(\mathbb C)$  is isomorphic to the adjoint $\mathfrak l$-module, where $\mathfrak l\simeq \mathfrak {gl}_{k_1|l_1}(\mathbb C)$ is the Lie superalgebra of $\mathcal L$. The rest of the proof is similar to the proof of Theorem \ref{teor osp}.$\Box$

\medskip

We have proved the following theorem.

\medskip

\t \label{teor general pisp} {\sl 
	Assume that $n=k_1+l_1$; the conditions (\ref{eq_condition no functions}) hold; $\mathfrak{v}(\mathbf{F}_{k'|l'})\simeq \mathfrak {pgl}_{k_1|l_1}(\mathbb
	C)$ and $k_1\geq 3$, $l_1\geq 2$. Then
	$\mathfrak v(\mathbf F_{k|l}^o)\simeq
	\pi\!\mathfrak {sp}_{n}(\mathbb C).\Box$}

\medskip

\bigskip

\bigskip

\bigskip

\noindent {\bf Acknowledgments.} The work was supported by the Universidade Federal de Minas Gerais, Brazil.

\bigskip

\bigskip

\noindent{\it Elizaveta Vishnyakova}

\noindent {Universidade Federal de Minas Gerais, Brazil}

\noindent{\emph{E-mail address:}
	\verb"VishnyakovaE@googlemail.com"}

\end{document}